\documentclass[12pt,reqno]{amsart}

\usepackage[english]{babel}
\usepackage[utf8]{inputenc}
\usepackage[T1]{fontenc}

\usepackage[margin=3cm]{geometry}
\usepackage{microtype}

\usepackage{mathtools,amssymb,mathrsfs}
\mathtoolsset{showonlyrefs}

\usepackage{graphicx}
\usepackage{subcaption}
\usepackage{booktabs}
\usepackage{caption}
\usepackage{tabularx}
\usepackage{array}
\usepackage{threeparttable}
\usepackage{adjustbox}
\usepackage{longtable}
\usepackage{pdflscape}
\usepackage{float}
\usepackage{placeins}

\usepackage{setspace}
\usepackage{multicol}
\usepackage{appendix}

\usepackage[authoryear,round]{natbib}

\usepackage[colorlinks=true,
            linkcolor=blue,
            citecolor=blue,
            urlcolor=blue]{hyperref}

\theoremstyle{plain}
\newtheorem{theorem}{Theorem}[section]

\newtheorem*{proposition*}{Proposition}
\newtheorem*{theorem*}{Theorem}

\newtheorem*{lemma*}{Lemma}

\newtheorem{assumption}{Assumption}[section]

\theoremstyle{definition}

\newcommand{\E}{\mathbb{E}}
\newcommand{\pr}{\mathbb{P}}

\newcommand{\cov}{\mathrm{Cov}}

\title{Self-Normalized Tests for Multistep Conditional Predictive Ability}

\author{Qitong Chen}
\address{School of Economics and Trade, Guangdong University of Foreign Studies, Guangzhou 510006, China}
\email{chenqitong@gdufs.edu.cn}

\author{Shuwen Lai}
\address{School of Statistics and Data Science, Nankai University, Tianjin 300071, China}
\email{laishuwen@mail.nankai.edu.cn}

\begin{document}

\begin{abstract}
This paper proposes self-normalized tests for multistep conditional predictive ability in forecast comparison. 
By normalizing the sample mean of the transformed loss differential using functionals of its cumulative sum (CUSUM) process, specifically an adjusted-range normalizer for scalars and a matrix normalizer for vectors, our approach avoids direct estimation of the long-run covariance matrix.
Consequently, it eliminates the need for the ad hoc bandwidth, kernel, and lag-truncation choices required by traditional methods. 
We establish the asymptotic theory for these statistics, deriving pivotal null limiting distributions and proving test consistency. 
Monte Carlo simulations show that the proposed tests effectively mitigate the finite-sample size distortions associated with traditional heteroskedasticity and autocorrelation consistent (HAC) methods, while retaining strong empirical power against conditional predictability alternatives.
\end{abstract}

\maketitle
\onehalfspacing

\section{Introduction}

Forecast comparison remains a central issue in time series econometrics. The traditional literature on forecast evaluation, such as \cite{DM95} and \cite{West96}, primarily focuses on testing the validity of specific forecasting models at the population level. In practice, however, it is often more relevant to compare forecasting methods rather than the underlying specifications of particular models. A forecasting method may encompass model-based, semiparametric, nonparametric, or algorithmic approaches; the primary object of interest is the sequence of forecasts generated from the information available at each forecast origin. The conditional predictive ability (CPA) framework proposed by \cite{GW06} is well suited to this perspective, as it tests whether relative forecast performance is predictable based on the information set available at the time the forecasts are formed.

While the CPA framework is versatile, evaluating multistep forecasts within this setup can introduce nontrivial statistical complications in finite samples. Analytically, overlapping forecast errors tend to give rise to an $\mathrm{MA}(\tau-1)$ serial correlation structure within the transformed loss differential. Consequently, the asymptotic null distribution of the test statistic depends on a complex nuisance parameter: a long-run variance in the scalar case, or a long-run covariance matrix in the vector-valued case. To address this nuisance parameter, standard CPA tests typically employ heteroskedasticity and autocorrelation consistent (HAC) estimation (\citealp{NW87,Andrews91}). However, HAC estimators necessitate the selection of bandwidths and kernel functions, and their convergence rates can be relatively slow. In finite samples, particularly when the forecast horizon $\tau$ is large or the underlying data exhibit high persistence, direct HAC estimation may yield ill-conditioned covariance matrices. This instability can lead to pronounced size distortions (i.e., over-rejection of the null hypothesis), which may constrain the reliability of traditional multistep CPA tests in macroeconomic and financial applications where evaluation samples are limited.

To circumvent the limitations associated with direct HAC estimation, we propose a self-normalized (SN) multistep CPA test. The SN approach, systematically developed by \cite{KV02, KV05} and \cite{Shao10}, replaces the explicit estimation of nuisance covariance matrices with a normalizer constructed from the recursive fluctuations (e.g., partial-sum processes) of the underlying data. While this paradigm has recently been extended to various complex settings, such as structural break testing (\citealp{Shao10, Hong2024}), high-dimensional time series inference (\citealp{WangShao2020, Wangetal2022}), and unconditional forecast evaluation (\citealp{Li2025}), accommodating the information-set-measurable test functions required for \textit{conditional} predictive ability presents an unresolved methodological challenge. Building upon the principle of self-normalization, we formulate our test by interacting the $\tau$-step-ahead loss differential, $\Delta L_{m,t+\tau}$, with an $\mathscr{F}_t$-measurable test function $h_t$, yielding a combined sequence $Z_{m,t+\tau} = h_t \Delta L_{m,t+\tau}$. The construction of the appropriate normalizer is governed by the dimensionality of the test function. For a scalar test function, we employ the adjusted range of the centered cumulative sum (CUSUM) process of $Z_{m,t+\tau}$. Conversely, for a vector-valued test function, for which extending the adjusted range is mathematically nontrivial due to cross-sectional dependence, we adopt a matrix CUSUM self-normalizer. A primary theoretical advantage of this methodology over traditional HAC procedures is that the limiting behavior of both the sample mean and the self-normalizer involves the same unknown long-run variance in the scalar case, or the same unknown long-run covariance matrix in the vector-valued case. Consequently, when the test statistic is constructed from these two objects, the nuisance variance or covariance parameter cancels out asymptotically. This mechanism yields a pivotal limiting distribution, thereby avoiding subjective bandwidth, kernel, and lag-truncation choices and improving finite-sample reliability.

This paper makes three main contributions to the literature on forecast evaluation. First, we develop a self-normalized conditional predictive ability (SN-CPA) testing framework applicable to both scalar and vector-valued $\mathscr{F}_t$-measurable test functions. By using functionals of the cumulative sum (CUSUM) process of the transformed loss differential, this method circumvents direct long-run covariance estimation and the subjective tuning parameters inherent in conventional HAC procedures. Second, we establish the asymptotic properties of the proposed test statistics. We show that the self-normalizers asymptotically account for the unknown long-run covariance matrix, yielding pivotal limiting null distributions. Furthermore, we prove test consistency against fixed alternatives and derive noncentral limits under local alternatives. Third, Monte Carlo simulations demonstrate that the SN-CPA tests systematically control the finite-sample size distortions commonly associated with HAC-based methods, especially over extended forecast horizons.

The rest of the paper is organized as follows. Section \ref{sec:setup} introduces the multistep CPA framework and the proposed statistics. Section \ref{sec:theory} presents the asymptotic theory. Section \ref{sec:simulation} reports the Monte Carlo evidence. Section \ref{sec:conclusion} concludes. Proofs are available upon request.

\section{Conditional predictive ability framework} \label{sec:setup}

This section introduces the multistep CPA framework and explains the choice of normalizers used in the main text. 

\subsection{Forecast comparison setup}
\label{subsec:forecast-setup}

Let $W_t=(Y_t,X_t^T)^T$ be the observed process, where $Y_t$ is the variable of interest and $X_t$ is a vector of predictor variables. 
Let $\mathscr F_t=\sigma(W_1,\ldots,W_t)$ denote the information set available at time $t$.
At forecast origin $t$, consider two competing forecasting methods for the
future value $Y_{t+\tau}$, where $\tau > 1$ is the forecast horizon. The
two forecasts are denoted by $\widehat f_t$ and $\widehat g_t$.
They may be generated by parametric, semiparametric, or nonparametric methods,
provided that both forecasts are $\mathscr F_t$-measurable.

We use a finite rolling-window estimation scheme for forecast evaluation. 
Let $m_{f,t}$ and $m_{g,t}$ denote the estimation-window lengths used by the two forecasting
methods at forecast origin $t$.
These window lengths may be method-specific constants or $\mathscr F_t$-measurable random integers determined by the forecasting methods.
Throughout the paper, we assume that they are uniformly bounded by a finite constant $\bar m$, namely, $\sup_t \max \{m_{f,t}, m_{g,t} \} \leq \bar{m} <\infty$.
Let $m$ denote the first forecast origin at which both methods can produce a $\tau$-step-ahead forecast.
The out-of-sample evaluation period is $t=m,\ldots,T-\tau$,  and hence $ n=T-\tau-m+1$.  
The asymptotics are taken with $m$ fixed and $n\to\infty$, so they are driven by the evaluation sample rather than by an increasing estimation window. 
To avoid excessive notation, we suppress the explicit dependence of the forecasts on $m_{f,t}$ and $m_{g,t}$.

Let $L(\cdot,\cdot)$ be a prespecified loss function, and define the $\tau$-step-ahead loss differential as $\Delta L_{m,t+\tau} = L(Y_{t+\tau},\widehat f_t) - L(Y_{t+\tau},\widehat g_t)$. The null hypothesis of equal conditional predictive ability (CPA) is formulated as 	
\begin{equation*}
H_0: \mathbb{E}[\Delta L_{m,t+\tau}\mid \mathscr{F}_t] = 0 \quad \text{almost surely for all} \ \ t.
\end{equation*}

Testing this conditional moment restriction directly against all $\mathscr{F}_t$-measurable alternatives is practically infeasible. Following \cite{GW06}, we introduce a $q\times 1$ vector of $\mathscr{F}_t$-measurable test functions, $h_t$, to project the conditional hypothesis into a tractable finite-dimensional space. Then, $H_0$ implies the necessary condition $\mathbb{E}[Z_{m,t+\tau}] = \mathbf{0}$, where $Z_{m,t+\tau} = h_t \Delta L_{m,t+\tau} \in \mathbb{R}^q$ is the transformed loss differential. Thus, the CPA framework operationalizes forecast conditional evaluation by testing whether the unconditional mean of $Z_{m,t+\tau}$ is zero.

The specification of $h_t$ is crucial, as it dictates the specific directions in which the test has power. 
When $h_t \equiv 1$, the projected restriction reduces to a conventional unconditional forecast comparison. Conversely, incorporating lagged loss differences or relevant macroeconomic state variables into $h_t$ directs the test's power toward detecting persistent or state-dependent predictability in relative forecast performance. In this paper, we focus on fixed and low-dimensional test functions. \footnote{While a high-dimensional $h_t$ encompasses a broader information set, it inevitably inflates the dimensionality of the normalizer, introducing estimation noise that  erodes finite-sample power. Exploring consistent conditional moment testing via infinite-dimensional spaces (e.g., \cite{Bie90}; \cite{SW98}) presents a different statistical trade-off and is beyond the scope of this paper.}

\subsection{Multistep dependence and choice of normalizers}
\label{subsec-normalizers}

\begin{table}[htbp]
\centering
\caption{Multistep regimes for self-normalized CPA testing}
\label{tab:multistep-cases}
\small
\renewcommand{\arraystretch}{1.15}
\begin{tabularx}{0.95\textwidth}{@{} l l X X @{}}
\toprule
Regime & Dependence & Nuisance object & Normalizer \\
\midrule
$q=1,\ \tau>1$ 
& scalar serial dependence 
& long-run variance 
& adjusted range \\

$q>1,\ \tau>1$ 
& vector serial dependence 
& long-run covariance matrix 
& matrix CUSUM self-normalizer \\
\bottomrule
\end{tabularx}
\end{table}

For $\tau>1$, overlapping forecast errors generally make $Z_{m,t+\tau}$ serially dependent under the null. 
Hence the nuisance object is a long-run variance when $q=1$, and a long-run covariance matrix when $q>1$. 
Standard CPA procedures handle this dependence through HAC estimation, which requires bandwidth, kernel, and lag-truncation choices. 
The purpose of the proposed normalizers is to avoid directly estimating this long-run variance or covariance matrix.

When $q=1$, $Z_{m,t+\tau}$ is scalar.  A conventional studentized statistic would require an estimate of the long-run variance. 
We instead use an adjusted-range normalizer based on the CUSUM process. 
This normalizer cancels out the unknown long-run variance without requiring HAC estimation. The use of the adjusted range is also in line with recent range-based self-normalized procedures for time-series testing, such as \cite{Hong2024}.

When $q>1$,  $Z_{m,t+\tau}$ is both vector-valued and serially dependent.
The nuisance object is then a full long-run covariance matrix.
HAC-based Wald procedures require estimating the precision matrix and are therefore sensitive to tuning choices and finite-sample estimation error.  
A componentwise range normalizer is not sufficient in this setting, because it does not account for the joint long-run dependence across coordinates. Moreover, the range-based multivariate extensions typically require ordering or decorrelation choices; see \cite{Hong2024}.  
We therefore use a matrix CUSUM self-normalizer, following the self-normalization principle of \cite{Shao10}, without introducing additional tuning-parameter choices. 
This construction absorbs the unknown long-run covariance matrix and avoids  HAC bandwidth, kernel, lag-truncation, and separate scale-estimation choices.

These considerations lead to the two multistep cases summarized in Table~\ref{tab:multi-cases}.

\section{Self-normalized statistics and asymptotic theory}
\label{sec:theory}

This section presents the self-normalized statistics and their asymptotic theory for multistep CPA testing. Throughout this section, $B$ denotes a standard one-dimensional Brownian motion and $\mathbb B(r)=B(r)-rB(1)$ denotes the associated Brownian bridge.
In the multivariate case, $B^q$ and $\mathbb{B}^q(r)=B^q(r)-rB^q(1)$ denote their $q$-dimensional counterparts.

We consider fixed and local alternatives. The fixed alternative is
\begin{equation*}
	H^{(q,\tau),F}_{A,h}: \mathbb E[Z_{m,t+\tau}\mid \mathscr F_t]=\mu_F \quad \text{almost surely for all} \ \  t, \qquad \mu_F\in\mathbb R^q,\ \mu_F\neq \mathbf{0},
\end{equation*}
and we write $V_{m,t+\tau}^{F}=Z_{m,t+\tau}-\mu_F$.
The local alternative is
\begin{equation*}
H^{(q,\tau),L}_{A,h}:  \mathbb E[Z_{m,t+\tau}\mid \mathscr F_t] =  n^{-1/2}\mu_L \quad \text{almost surely for all} \ \ t, \qquad \mu_L\in\mathbb R^q,\ \mu_L\neq \mathbf{0},
\end{equation*}
and we write $V_{m,t+\tau}^{L}=Z_{m,t+\tau}-n^{-1/2}\mu_L$.  When $q=1$, $\mu_F$ and $\mu_L$ are understood as nonzero real constants.

The fixed and local alternatives describe the power behavior of the proposed statistics.
They are stronger than the global alternatives in \cite{GW06}, but have a direct CPA interpretation: $h_t$ contains information that predicts relative forecast performance. 
 A fixed drift makes the statistic diverge, whereas an $n^{-1/2}$-local drift enters the limiting self-normalized Brownian functional as a noncentrality parameter. 
 In both cases, the CUSUM normalizer removes these drifts from the denominator but not from the numerator.

\noindent \textbf{\textit{Remark} (State-dependent alternatives).} 
	\label{rem:state-dependent-alt}
	The fixed and local alternatives above are stated with constant conditional drifts for notational simplicity. 
	The same conclusions continue to hold under state-dependent conditional drifts, provided that the average drift is nonzero under fixed alternatives and the average local drift coefficient is nonzero under local alternatives, while the centered drift components do not affect the centered CUSUM terms used in the normalizers below.
	
Let $g_{m,t+\tau}:=\E[Z_{m,t+\tau}\mid \mathscr F_t]$ and $\overline g_{m,n,\tau}:=n^{-1}\sum_{t=m}^{m+n-1}g_{m,t+\tau}$. Consider a fixed alternative under which $\overline g_{m,n,\tau} \overset{p}{\longrightarrow} \mu_F$ for some $\mu_F\in\mathbb R^q$ with $\mu_F\neq \mathbf{0}$. Define $V_{m,t+\tau}^{F}:= Z_{m,t+\tau}-g_{m,t+\tau}$. Assume that the FCLT continues to hold when $Z_{m,t+\tau}$ is replaced by $V_{m,t+\tau}^{F}$. Suppose further that the centered drift component satisfies $\sup_{0\le r\le 1}\left\|n^{-1/2}\sum_{t=m}^{m+\lfloor nr\rfloor-1}\left(g_{m,t+\tau}-\overline g_{m,n,\tau}\right)\right\|=O_p(1)$. Then the centered CUSUM term based on $Z_{m,t+\tau}$ remains of stochastic order $O_p(1)$. Indeed, it can be decomposed into the centered CUSUM term based on $V_{m,t+\tau}^{F}$, which is $O_p(1)$ by the FCLT, and the centered drift component above. On the other hand, $\overline Z_{m,n,\tau}=\overline g_{m,n,\tau}+n^{-1}\sum_{t=m}^{m+n-1}V_{m,t+\tau}^{F}=\mu_F+o_p(1)$. Thus, the fixed-alternative consistency conclusion stated below for constant drifts continues to hold under this state-dependent fixed alternative.

For local alternatives, let $g_{m,t+\tau}=n^{-1/2}a_{m,t+\tau}$ and $\overline a_{m,n,\tau}:=n^{-1}\sum_{t=m}^{m+n-1}a_{m,t+\tau}\overset{p}{\longrightarrow}\mu_L$,
for some $\mu_L\in\mathbb R^q$ with $\mu_L\neq\mathbf{0}$.
Define $V_{m,t+\tau}^{L}:=Z_{m,t+\tau}-n^{-1/2}a_{m,t+\tau}$. Assume analogously that the FCLT continues to hold with $Z_{m,t+\tau}$ replaced by $V_{m,t+\tau}^{L}$, and that $\sup_{0\le r\le 1}\left\|n^{-1/2}\sum_{t=m}^{m+\lfloor nr\rfloor-1}\left(a_{m,t+\tau}-\overline a_{m,n,\tau}\right)\right\|=O_p(1)$. Then the same decomposition argument shows that the local drift is removed from the centered CUSUM term but remains in the numerator. Hence, the local noncentral limiting conclusions stated below for constant local drifts continue to hold under this state-dependent local alternative, with the noncentrality parameter determined by the limit $\mu_L$.

\subsection{Scalar multistep case: $q=1$}
\label{subsec:scalar-multi-theory}
Define the CUSUM process
\begin{equation*}
	 T^h_{m,n,\tau}(r) = \frac{1}{\sqrt n} \sum_{t=m}^{m+\lfloor nr\rfloor-1}
        \bigl(Z_{m,t+\tau}-\overline Z_{m,n,\tau}\bigr), \qquad 0\leq r\leq 1,
\end{equation*}
where  $ \overline Z_{m,n,\tau} = n^{-1}\sum_{t=m}^{T-\tau} Z_{m,t+\tau}$ and $\lfloor x \rfloor =\max \{z\leq x \colon z\in\mathbb{Z} \}$.
The adjusted range is $R^h_{m,n,\tau}= \sup_{0\leq r\leq 1}T^h_{m,n,\tau}(r)- \inf_{0\leq r\leq 1}T^h_{m,n,\tau}(r)$, and the scalar multistep CPA test statistic is 
\begin{equation*}
        Q^h_{m,n,\tau}=\frac{n\overline Z_{m,n,\tau}^{\,2}}{(R^h_{m,n,\tau})^2}, \qquad \text{whenever $R^h_{m,n,\tau}>0$.}
\end{equation*}

To show the limiting distribution under $H_0$, we impose the following regularity conditions.
\begin{assumption} \label{ass:scalar-multi}
(1)	The sequence $\{Z_{m, t+\tau} \colon t\geq m \}$ is strongly $\alpha$-mixing with coefficients $\alpha(k)$, and there exists $\eta>0$ such that  $\sum_{k=1}^{\infty} \alpha(k)^{\eta/(2+\eta)}<\infty$.
(2) For the same $\eta>0$,  $ \sup_{t\geq m}\E|Z_{m, t+\tau}|^{2+\eta}<\infty$.
(3) There exists $\sigma_{\tau}^2 \in (0,\infty)$ such that, for every $r,s \in[0,1]$, $n^{-1}\cov\left( \sum_{t=m}^{m+\lfloor nr\rfloor-1} Z_{m, t+\tau}, \sum_{t=m}^{m+\lfloor ns\rfloor-1} Z_{m, t+\tau}\right) \longrightarrow  (r\wedge s)\sigma_{\tau}^2$.
\end{assumption}

Assumption~\ref{ass:scalar-multi} is a high-level regularity condition imposed directly on  the sequence $\{Z_{m,t+\tau}=h_t\Delta L_{m,t+\tau},\ t\geq m\}$, where $m$ is the starting point of the forecast evaluation period. These conditions can be verified from primitive assumptions on the underlying process $W_t$, the forecasting methods, the loss function, and the test function $h_t$, as in \cite{GW06}; for details, see the Supplementary Material. Moreover, Assumption~\ref{ass:scalar-multi} is consistent with standard strong-mixing invariance principles for the dependent sequence; see \cite{Herrndorf84,Her85} and \cite{deJ00}. 
Theorem \ref{thm:scalar-multi-null} below establishes pivotal null limits.

\begin{theorem}\label{thm:scalar-multi-null}
Suppose Assumption~\ref{ass:scalar-multi} holds. Then, under \(H_0\), 
\begin{equation*}
Q^h_{m,n,\tau} \Rightarrow  \frac{B(1)^2} {\left( \sup_{0\leq r\leq 1}\mathbb B(r)-\inf_{0\leq r\leq 1}\mathbb B(r) \right)^2}.	
\end{equation*}
\end{theorem}

Under Assumption~\ref{ass:scalar-multi} and \(H_0\), the partial sum process satisfies
\begin{equation}
	n^{-1/2}\sum_{t=m}^{m+\lfloor n\cdot\rfloor-1}Z_{m,t+\tau} \Rightarrow\sigma_\tau B(\cdot)\ \text{in } D([0,1],\mathbb R).
\end{equation} 
Hence, by the continuous mapping theorem (CMT),  $T^h_{m,n,\tau}(\cdot)\Rightarrow \sigma_\tau \mathbb B(\cdot)$ and $R^h_{m,n,\tau}\Rightarrow\sigma_\tau\left(\sup_{0\le r\le 1}\mathbb B(r)-\inf_{0\le r\le 1}\mathbb B(r)\right)$.
Since the limiting range is strictly positive almost surely, $\pr(R^h_{m,n,\tau}>0)\to 1$.
Hence \(Q^h_{m,n,\tau}\) is well defined with probability approaching one.

As demonstrated in Theorem \ref{thm:scalar-multi-null}, the unknown long-run standard deviation $\sigma_\tau$ cancels out asymptotically, yielding a pivotal null limit. Unlike the standard HAC-based procedures in \cite{GW06}, this self-normalized approach bypasses explicit variance estimation and its associated tuning sensitivities. The adjusted range is motivated by recent range-based self-normalized procedures, such as \cite{Hong2024}, and we adapt it to multistep CPA testing with persistent loss differentials.

To evaluate the power properties of the proposed test, Theorem \ref{thm:scalar-multi-alt} below establishes its consistency against fixed alternatives and derives the noncentral limit under local alternatives. 

\begin{theorem}
\label{thm:scalar-multi-alt}
Suppose the statistic $Q^h_{m,n,\tau}$ is defined as above.
\begin{enumerate}
\item[(i)] Under $H^{(1,\tau),F}_{A,h}$, if Assumption~\ref{ass:scalar-multi} holds with $Z_{m,t+\tau}$ replaced by $V_{m,t+\tau}^{F}$, then $ Q^h_{m,n,\tau}\xrightarrow{p}\infty$.
\item[(ii)] Under \(H^{(1,\tau),L}_{A,h}\), if Assumption~\ref{ass:scalar-multi} holds with $Z_{m,t+\tau}$ replaced by $V_{m,t+\tau}^{L}$, then, with $J_\tau=\mu_L/\sigma_\tau$, 
\begin{equation*}
	Q^h_{m,n,\tau}\Rightarrow \frac{(B(1)+J_\tau)^2}{\left(\sup_{0\leq r\leq 1}\mathbb B(r) -\inf_{0\leq r\leq 1}\mathbb B(r) \right)^2}.
\end{equation*}
\end{enumerate}
\end{theorem}

 Unlike \cite{GW06}, which focuses mainly on fixed alternatives, Theorem \ref{thm:scalar-multi-alt} also derives the noncentral limit under local alternatives and thereby  characterizes local asymptotic power.

\subsection{Multivariate multistep case: $q>1$} 
\label{subsec:vector-multi-theory}

In this subsection, we use a matrix CUSUM self-normalizer, following the self-normalization principle of \cite{Shao10}, rather than range normalization because of the joint long-run dependence across coordinates.  

Define 
 \begin{equation*}
	T_{m,n,\tau}^{h,q}(r) =\frac{1}{\sqrt{n}} \sum_{t=m}^{m+\lfloor nr\rfloor-1} \bigl(Z_{m, t+\tau}-\overline Z_{m,n,\tau}\bigr), \qquad 0\leq r\leq 1,
\end{equation*}
where $\overline Z_{m,n,\tau}= n^{-1} \sum_{t=m}^{T-\tau}Z_{m, t+\tau}$.
The matrix self-normalizer is 
\begin{equation*}
	U_{m,n,\tau}= n^{-1} \sum_{k=1}^{n} T^{h,q}_{m,n,\tau}\!\left(k/n\right) T^{h,q}_{m,n,\tau}\!\left(k/n \right)^T.
\end{equation*}
The multivariate multistep CPA test statistic is 
\begin{equation*}
        Q^{h,q}_{m,n,\tau}=
        n\overline Z_{m,n,\tau}^T
        U_{m,n,\tau}^{-1}
        \overline Z_{m,n,\tau}, \qquad \text{whenever \(U_{m,n,\tau}\) is nonsingular.}
\end{equation*} 

\begin{assumption} \label{ass:vector-multi}
(1) The sequence $\{Z_{m, t+\tau} \colon t\geq m \}$ is strongly $\alpha$-mixing with coefficients $\alpha(k)$, and there exists $\eta>0$ such that  $\sum_{k=1}^{\infty} \alpha(k)^{\eta/(2+\eta)}<\infty$.
(2) For the same $\eta>0$, $\sup_{t\geq m}  \E\|Z_{m,t+\tau}\|^{2+\eta}<\infty$.
(3) There exists a positive definite $q\times q$ matrix $\Sigma_\tau$ such that, for every \(r,s\in[0,1]\), $n^{-1} \cov \left( \sum_{t=m}^{m+\lfloor nr\rfloor-1}Z_{m,t+\tau}, \sum_{t=m}^{m+\lfloor ns\rfloor-1}Z_{m,t+\tau}  \right)  \longrightarrow (r\wedge s)\Sigma_\tau $.
\end{assumption}        
Assumption~\ref{ass:vector-multi} is a high-level regularity condition imposed directly on $Z_{m,t+\tau}$. Since $m$ and $\tau$ are fixed, this process is treated as an ordinary $\mathbb{R}^q$-valued strongly mixing sequence indexed by $t\geq m$. 
The mixing, moment, and covariance conditions are sufficient to apply the strong-mixing FCLT; see \cite{Herrndorf84,Her85} and \cite[Chapter~29]{Davidson94}. Theorem~\ref{thm:vector-multi-null} establishes the null limit of $Q^{h,q}_{m,n,\tau}$.

\begin{theorem} \label{thm:vector-multi-null}
Suppose Assumption~\ref{ass:vector-multi} holds. Then, under $H_0$,  
\begin{equation*}
	Q^{h,q}_{m,n,\tau}\Rightarrow  B^q(1)^T \left(  \int_0^1 \mathbb B^q(r)\mathbb B^q(r)^T dr  \right)^{-1} B^q(1).
\end{equation*}
\end{theorem}

The multivariate case is analogous to the scalar case,
with the range functional replaced by the matrix CUSUM self-normalizer.
Under Assumption~\ref{ass:vector-multi} and \(H_0\), 
\begin{equation}
n^{-1/2}\sum_{t=m}^{m+\lfloor n\cdot\rfloor-1}Z_{m,t+\tau} \Rightarrow \Sigma_\tau^{1/2}B^q(\cdot) \ \text{in }D([0,1],\mathbb R^q).	
\end{equation}
Hence, by the CMT and the Riemann sum representation of $U_{m,n,\tau}$, $T^{h,q}_{m,n,\tau}(\cdot)\Rightarrow\Sigma_\tau^{1/2}\mathbb B^q(\cdot)$ and $U_{m,n,\tau}\Rightarrow\Sigma_\tau^{1/2}\left(\int_0^1\mathbb B^q(r)\mathbb B^q(r)^\top\,dr\right)\Sigma_\tau^{1/2}$.
Since the limiting matrix is nonsingular almost surely, \(U_{m,n,\tau}\) is nonsingular with probability approaching one, and hence \(Q^{h,q}_{m,n,\tau}\) is well defined with probability approaching one.

Unlike \cite{GW06}, this pivotal limit in Theorem \ref{thm:vector-multi-null} is achieved without the explicit estimation of the long-run precision matrix. The matrix CUSUM self-normalizer asymptotically absorbs the nuisance covariance matrix $\Sigma_\tau$, thereby extending the tuning-free principle of \cite{Shao10} to multistep forecast comparisons. 

Having established the pivotal null limit for the multivariate test, Theorem \ref{thm:vector-multi-alt} characterizes its asymptotic power against both fixed and local alternatives.

\begin{theorem}
\label{thm:vector-multi-alt}
Suppose the statistic $Q^{h,q}_{m,n,\tau}$ is defined as above.
\begin{enumerate}
\item[(i)] Under $H^{(q,\tau),F}_{A,h}$, if Assumption~\ref{ass:vector-multi} holds with $Z_{m,t+\tau}$ replaced by $V_{m,t+\tau}^{F}$, then $Q^{h,q}_{m,n,\tau}\xrightarrow{p}\infty $.
\item[(ii)] Under $H^{(q,\tau),L}_{A,h}$, if Assumption~\ref{ass:vector-multi} holds with $Z_{m,t+\tau}$ replaced by $V_{m,t+\tau}^{L}$, then, with $ J_\tau=\Sigma_\tau^{-1/2}\mu_L$, we have 
\begin{equation*}
	 Q^{h,q}_{m,n,\tau}  \Rightarrow (B^q(1)+J_\tau)^T\left(\int_0^1  \mathbb B^q(r)\mathbb B^q(r)^T dr  \right)^{-1}(B^q(1)+J_\tau).
\end{equation*}
\end{enumerate}
\end{theorem}

The matrix self-normalizer asymptotically absorbs the nuisance covariance $\Sigma_\tau$ under the null hypothesis and implicitly standardizes the local drift under the alternative. Consequently, the noncentrality vector $J_\tau$ represents the local drift standardized by the long-run covariance, while preserving local asymptotic power without requiring explicit HAC estimation.

\section{Monte Carlo simulation}\label{sec:simulation}
This section reports Monte Carlo simulation results for the finite-sample performance of the proposed self-normalized CPA tests. For notational simplicity in the simulation and empirical sections, we denote the scalar statistic $Q^h_{m,n,\tau}$ as $Q_1$, and the multivariate statistic $Q^{h,q}_{m,n,\tau}$ as $Q_q$ (specifically, $Q_2$ when $q=2$).

\subsection{Data-generating processes (DGPs)}

We consider two DGPs. For each DGP, size and power are evaluated under the same structural framework by varying a single control parameter. 
Empirical size is computed when this parameter makes the conditional mean of the loss
differential zero, and power is computed when it makes the conditional mean nonzero.

\textbf{DGP 1: Conditional relative performance with overlapping errors.}
In the spirit of \cite{GW06}, we generate the predictor \(x_t\), used in the test function \(h_t=(1,x_t)^T\) for the multivariate test and \(h_t=x_t\) for the scalar test, as a persistent AR(1) process: $x_t=\rho x_{t-1}+u_t$, where $u_t\sim \text{i.i.d.}\ N(0,1)$ and $\rho \in \{0.2, 0.5, 0.8\}$ controls the persistence level. The sequence \(x_t\) is standardized to maintain a unit unconditional variance. 

The transformed loss differential is generated by $\Delta L_{m,t+\tau}=\delta x_t+\varepsilon_{t+\tau}$, where the forecast error $\varepsilon_{t+\tau}$ is modeled as an $\mathrm{MA}(\tau-1)$ process to explicitly capture the overlap in multistep forecasts:
$
\varepsilon_{t+\tau}
=
c \sum_{j=0}^{\tau-1}\theta_j v_{t+\tau-j},
$
with \(\theta_0=1\) and \(\theta_j=0.5\) for \(j>0\). The scaling constant \(c\) ensures that the unconditional variance of \(\varepsilon_{t+\tau}\) is normalized to one.

When \(\delta=0\), we have $\mathbb E[\Delta L_{m,t+\tau}\mid \mathscr F_t]=0$, so the design evaluates empirical size under the exact null hypothesis. When \(\delta>0\), $ \mathbb E[\Delta L_{m,t+\tau}\mid \mathscr F_t]=\delta x_t$. Importantly, because $\mathbb E[x_t]=0$, the unconditional mean $\mathbb E[\Delta L_{m,t+\tau}]$ remains identically zero regardless of \(\delta\). In this paper, we set $\delta=0.2$. Thus, relative performance is conditionally predictable but not unconditionally biased, thereby isolating the power against conditional alternatives while demonstrating the failure of unconditional benchmark tests.

\textbf{DGP 2: State-dependent performance.}
To evaluate the tests' sensitivity to regime shifts and business cycle fluctuations, we consider a state-dependent framework. The state variable $S_t \in \{0,1\}$ is assumed to be in the information set $\mathscr{F}_t$ and is drawn from a Bernoulli distribution with $P(S_t=1)=p$, where the state probability $p \in \{0.2, 0.5, 0.8\}$. The test functions are defined as $h_t=(1,S_t)^T$ for the multivariate test and $h_t=S_t$ for the scalar test.

The transformed loss differential is generated by: $\Delta L_{m,t+\tau} = d(S_t-p) + \varepsilon_{t+\tau}$, where the control parameter $d$ dictates the magnitude of state-dependent predictability, and $\varepsilon_{t+\tau}$ follows the same overlapping $\mathrm{MA}(\tau-1)$ structure as defined above.

When $d=0$, the conditional mean satisfies $\mathbb E[\Delta L_{m,t+\tau}\mid \mathscr F_t]=0$, establishing the exact null hypothesis for empirical size evaluation. When $d>0$, the relative forecast performance varies with the realized state $S_t$, so that the conditional expectation becomes $\mathbb E[\Delta L_{m,t+\tau}\mid \mathscr F_t]=d(S_t-p)$. 

Crucially, by construction, the unconditional expectation of the state indicator is $\mathbb E[S_t]=p$. Consequently, the unconditional mean of the loss differential remains identically zero ($\mathbb E[\Delta L_{m,t+\tau}] = 0$) for all values of $d$. This deliberately asymmetric design ensures that while the competing models exhibit strong conditional predictability across different states, they appear equally accurate unconditionally.

\subsection{Implementation and test statistics}

For both DGPs, we consider evaluation sample sizes \(n \in \{50, 100, 150, 200, 250, 300, 350, 400\}\) and perform \(5{,}000\) Monte Carlo replications. Empirical rejection frequencies are evaluated at the 1\%, 5\%, and 10\% nominal significance levels. To capture the conditional information, we employ the scalar test function \(h_t = x_t\) (for DGP 1) and \(h_t = S_t\) (for DGP 2) in the one-dimensional tests. For the multidimensional tests, we naturally expand the information set to include an intercept, utilizing the vector-valued test functions \(h_t = (1, x_t)^T\) and \(h_t = (1, S_t)^T\), respectively.

We systematically compare five test statistics to benchmark the finite-sample performance: The proposed vector-valued SN-CPA statistic (\(Q_{2}\)); the proposed scalar SN-CPA statistic (\(Q_{1}\)); the unconditional self-normalized DM statistic of \cite{Li2025} (\(T_{\mathrm{SN}}\)); the traditional HAC-based conditional predictive ability Wald statistic of \cite{GW06} (\(T_{\mathrm{GW}}\)); and the traditional HAC-based unconditional Diebold-Mariano statistic of \cite{DM95} (\(T_{\mathrm{DM}}\)).

Note that all HAC covariance matrices are estimated using the \texttt{NeweyWest} function from the \texttt{sandwich} package in R, maintaining all default settings. Because the asymptotic null distributions of our proposed SN-CPA statistics (\(Q_{1}\) and \(Q_{2}\)) are non-standard functionals of Brownian motions and Brownian bridges, we simulate their critical values by Monte Carlo methods. Following the theoretical representations in Section \ref{sec:theory}, we approximate the continuous-time Wiener processes using discrete Gaussian random walks on a fine grid with \(N=200,000\) steps and $M = 10{,}000$ independent replications. The critical values are reported in Table \ref{tab:critical_values_multistep}. Moreover, the critical values for the benchmark statistics (\(T_{\mathrm{SN}}\), \(T_{\mathrm{GW}}\), and \(T_{\mathrm{DM}}\)) are drawn from their well-established asymptotic distributions or existing literature.

\begin{table}[h]
	\centering
	\caption{Quantiles of the limiting null distribution ($\tau \ge 2$)}
	\label{tab:critical_values_multistep}
	\renewcommand{\arraystretch}{1.2}
	\resizebox{\textwidth}{!}{
		\begin{tabular}{l ccccccccc}
			\toprule
			& \multicolumn{9}{c}{Quantiles} \\
			\cmidrule(lr){2-10}
			Dimension ($q$) & 1\% & 5\% & 10\% & 25\% & 50\% & 75\% & 90\% & 95\% & 99\% \\
			\midrule
			$q=1$  & 0.000	 &0.003	&0.010	&0.066	&0.294	&0.934	&2.011	&3.071	&6.121\\			
			$q=2$          & 0.215 & 0.979 & 1.999 & 5.590 & 15.338 & 36.449 & 72.375 & 103.114 & 196.967 \\
			$q=3$          & 1.367 & 4.150 & 6.978 & 15.518 & 34.994 & 71.464 & 128.131 & 177.017 & 307.613 \\
			$q=4$          & 3.908 & 9.825 & 15.178 & 31.212 & 61.731 & 116.810 & 190.662 & 258.421 & 428.807 \\
			$q=5$          & 8.930 & 19.056 & 28.625 & 52.214 & 96.632 & 170.682 & 270.727 & 354.125 & 555.283 \\
			\bottomrule
			\multicolumn{10}{p{0.98\textwidth}}{\footnotesize \textit{Notes:} This table reports the asymptotic quantiles of the proposed self-normalized CPA test under multistep forecast horizons ($\tau \ge 2$, including $\tau=2$ and $\tau=3$). The vector-valued cases ($q \ge 2$) are based on the matrix CUSUM normalizer. Critical values are approximated using discrete Gaussian random walks with $N=200{,}000$ steps and $M=10{,}000$ Monte Carlo replications.}
		\end{tabular}
	}
\end{table}

\subsection{Empirical size and power}

Tables \ref{tab:dgp1-size-tau2}, \ref{tab:dgp1-size-tau3}, and Figure \ref{fig:dgp1-power} report the empirical size and power ($\delta=0.2$) for DGP 1. The results for DGP 2 yield highly similar conclusions; to conserve space, we focus our discussion on DGP 1 in the main text and relegate the comprehensive tables and figures for DGP 2 to the Supplementary Material. By comparing the proposed statistics with traditional benchmarks, we uncover the following key empirical findings:

First, traditional HAC-based tests ($T_{\mathrm{GW}}$ and $T_{\mathrm{DM}}$) tend to over-reject the null hypothesis in finite samples, particularly as the forecast horizon $\tau$ and data persistence $\rho$ increase. In contrast, the proposed statistic $Q_{2}$ effectively mitigates these finite-sample size distortions, maintaining accurate size control near the nominal level across all considered sample sizes and horizons. Second, the unconditional tests ($T_{\mathrm{SN}}$ and $T_{\mathrm{DM}}$) display flat power curves hovering around 5\%, which naturally aligns with their theoretical design focusing on average performance rather than dynamic, information-set-based predictability. Meanwhile, the proposed $Q_{1}$ and $Q_{2}$ tests exhibit robust, monotonically increasing power that converges to one as the sample size $n$ grows. Third, although the HAC-based $T_{\mathrm{GW}}$ test occasionally reports numerically higher raw power than $Q_{2}$, this apparent advantage should be interpreted with caution, as it is partially driven by its empirical size inflation under the null. When accounting for size control, $Q_{2}$ provides a more balanced and reliable measure of true statistical sensitivity. Finally, while $Q_{1}$ can detect specific directional predictability, it is empirically less stable than its multidimensional counterpart. By augmenting the test function with an intercept ($h_t = (1, x_t)^T$), $Q_{2}$ preserves the full joint conditional moment structure. This prevents information loss, yielding consistently tighter size control and a more comprehensive evaluation of predictive ability.

These empirical findings align with our established asymptotic theory. The finite-sample challenges associated with HAC methods illustrate the inherent difficulty of explicitly estimating the long-run covariance matrix induced by overlapping $MA(\tau-1)$ forecast errors. Conversely, the accurate size control and robust power of $Q_{2}$ validate our primary theoretical claim: the matrix CUSUM self-normalizer absorbs the complex serial and cross-sectional dependence. By entirely bypassing ad hoc bandwidth selections, our proposed methodology offers a reliable and tuning-free alternative for finite-sample forecast evaluation.

\begin{table}[!htbp]
	\centering
	\begin{threeparttable}
		\caption{Empirical size under DGP 1 for \(\tau=2\)}
		\label{tab:dgp1-size-tau2}

		\fontsize{7.5}{8.3}\selectfont
		\setlength{\tabcolsep}{5.0pt}
		\renewcommand{\arraystretch}{1.02}

		\begin{tabular*}{\textwidth}{@{\extracolsep{\fill}}cc ccccc@{}}
			\toprule
			\(\rho\) & \(n\)
			& \(Q_{2}\) & \(Q_{1}\) & \(T_{\mathrm{SN}}\)
			& \(T_{\mathrm{GW}}\) & \(T_{\mathrm{DM}}\) \\
			\midrule
			0.2 & 50
			& .110/.061/.012 & .158/.087/.021 & .118/.065/.015
			& .212/.133/.047 & .199/.127/.048 \\
			0.2 & 100
			& .106/.055/.011 & .142/.075/.015 & .111/.056/.014
			& .183/.114/.039 & .190/.117/.041 \\
			0.2 & 150
			& .097/.050/.009 & .130/.073/.016 & .104/.055/.012
			& .170/.104/.033 & .169/.104/.039 \\
			0.2 & 200
			& .099/.053/.010 & .135/.072/.015 & .108/.058/.012
			& .177/.105/.033 & .178/.110/.036 \\
			0.2 & 250
			& .101/.052/.010 & .127/.066/.015 & .106/.057/.010
			& .173/.105/.030 & .175/.106/.036 \\
			0.2 & 300
			& .100/.055/.011 & .126/.065/.012 & .109/.054/.013
			& .162/.097/.030 & .177/.105/.035 \\
			0.2 & 350
			& .101/.049/.009 & .127/.067/.013 & .102/.048/.009
			& .158/.089/.026 & .162/.093/.029 \\
			0.2 & 400
			& .105/.056/.010 & .121/.066/.015 & .105/.053/.011
			& .158/.091/.027 & .172/.101/.034 \\
			\addlinespace[0.35em]
			0.5 & 50
			& .117/.063/.014 & .184/.105/.026 & .112/.062/.017
			& .235/.150/.055 & .195/.126/.049 \\
			0.5 & 100
			& .105/.055/.014 & .147/.075/.020 & .118/.065/.017
			& .201/.126/.040 & .190/.118/.046 \\
			0.5 & 150
			& .104/.052/.012 & .144/.077/.017 & .103/.055/.011
			& .187/.113/.034 & .175/.109/.036 \\
			0.5 & 200
			& .113/.058/.012 & .150/.077/.019 & .108/.063/.011
			& .193/.115/.035 & .180/.114/.033 \\
			0.5 & 250
			& .106/.054/.010 & .145/.071/.016 & .102/.051/.011
			& .189/.112/.031 & .173/.107/.035 \\
			0.5 & 300
			& .095/.052/.013 & .131/.067/.015 & .097/.049/.012
			& .171/.103/.032 & .167/.102/.033 \\
			0.5 & 350
			& .101/.057/.011 & .132/.063/.015 & .104/.052/.012
			& .180/.103/.032 & .172/.101/.029 \\
			0.5 & 400
			& .096/.044/.007 & .130/.066/.012 & .105/.049/.008
			& .168/.098/.028 & .162/.098/.030 \\
			\addlinespace[0.35em]
			0.8 & 50
			& .110/.058/.013 & .190/.101/.024 & .116/.063/.013
			& .245/.163/.056 & .202/.126/.049 \\
			0.8 & 100
			& .105/.052/.010 & .165/.090/.017 & .107/.051/.011
			& .212/.133/.042 & .180/.115/.040 \\
			0.8 & 150
			& .102/.051/.010 & .150/.076/.020 & .100/.052/.012
			& .196/.124/.043 & .167/.103/.038 \\
			0.8 & 200
			& .099/.048/.008 & .133/.071/.015 & .110/.059/.012
			& .192/.121/.036 & .176/.107/.036 \\
			0.8 & 250
			& .100/.050/.010 & .146/.072/.013 & .104/.056/.010
			& .201/.122/.041 & .178/.109/.037 \\
			0.8 & 300
			& .092/.048/.010 & .139/.070/.015 & .101/.053/.012
			& .198/.118/.034 & .168/.106/.035 \\
			0.8 & 350
			& .105/.051/.009 & .133/.065/.012 & .105/.060/.013
			& .204/.120/.039 & .177/.108/.034 \\
			0.8 & 400
			& .100/.050/.009 & .133/.064/.012 & .102/.053/.010
			& .196/.116/.035 & .170/.103/.030 \\
			\bottomrule
		\end{tabular*}

		\vspace{0.35em}
		\begin{tablenotes}[flushleft]
			\scriptsize
			\item \emph{Notes.} The table reports empirical rejection frequencies under the null with \(B=5{,}000\) replications. Each entry is reported as \(10\%/5\%/1\%\), corresponding to the nominal significance level \(\alpha\). \(Q_{2}\) denotes the proposed vector-valued SN-CPA statistic, \(Q_{1}\) denotes the proposed scalar SN-CPA statistic, \(T_{\mathrm{SN}}\) denotes the unconditional self-normalized DM statistic, \(T_{\mathrm{GW}}\) denotes the HAC-based CPA statistic, and \(T_{\mathrm{DM}}\) denotes the HAC-based DM statistic.
		\end{tablenotes}
	\end{threeparttable}
\end{table}

\begin{table}[!htbp]
	\centering
	\begin{threeparttable}
		\caption{Empirical size under DGP 1 for \(\tau=3\)}
		\label{tab:dgp1-size-tau3}

		\fontsize{7.5}{8.3}\selectfont
		\setlength{\tabcolsep}{5.0pt}
		\renewcommand{\arraystretch}{1.02}

		\begin{tabular*}{\textwidth}{@{\extracolsep{\fill}}cc ccccc@{}}
			\toprule
			\(\rho\) & \(n\)
			& \(Q_{2}\) & \(Q_{1}\) & \(T_{\mathrm{SN}}\)
			& \(T_{\mathrm{GW}}\) & \(T_{\mathrm{DM}}\) \\
			\midrule
			0.2 & 50
			& .106/.054/.011 & .166/.091/.023 & .117/.061/.013
			& .185/.112/.036 & .168/.099/.035 \\
			0.2 & 100
			& .102/.053/.010 & .138/.079/.019 & .103/.053/.011
			& .161/.091/.022 & .157/.092/.025 \\
			0.2 & 150
			& .096/.048/.009 & .127/.068/.016 & .095/.051/.012
			& .149/.086/.024 & .145/.087/.025 \\
			0.2 & 200
			& .100/.051/.010 & .131/.065/.015 & .102/.053/.010
			& .160/.090/.024 & .151/.092/.025 \\
			0.2 & 250
			& .099/.053/.011 & .127/.067/.014 & .104/.051/.010
			& .143/.081/.021 & .143/.084/.023 \\
			0.2 & 300
			& .109/.056/.012 & .128/.067/.015 & .102/.050/.011
			& .155/.082/.023 & .158/.090/.026 \\
			0.2 & 350
			& .099/.056/.010 & .121/.065/.015 & .103/.047/.010
			& .149/.084/.021 & .150/.088/.025 \\
			0.2 & 400
			& .108/.060/.013 & .123/.067/.012 & .111/.063/.014
			& .161/.087/.024 & .163/.097/.025 \\
			\addlinespace[0.35em]
			0.5 & 50
			& .107/.059/.011 & .178/.097/.026 & .112/.059/.014
			& .193/.118/.039 & .166/.100/.034 \\
			0.5 & 100
			& .105/.055/.012 & .154/.085/.021 & .106/.054/.012
			& .179/.103/.031 & .157/.095/.029 \\
			0.5 & 150
			& .098/.049/.012 & .138/.070/.016 & .103/.051/.010
			& .165/.090/.028 & .153/.088/.027 \\
			0.5 & 200
			& .103/.052/.010 & .143/.077/.017 & .101/.055/.011
			& .169/.101/.028 & .153/.092/.027 \\
			0.5 & 250
			& .102/.053/.011 & .130/.070/.014 & .111/.054/.010
			& .159/.095/.029 & .151/.083/.029 \\
			0.5 & 300
			& .101/.048/.010 & .130/.067/.013 & .105/.053/.011
			& .166/.093/.026 & .150/.090/.027 \\
			0.5 & 350
			& .096/.050/.008 & .121/.060/.012 & .099/.052/.012
			& .156/.090/.021 & .151/.091/.023 \\
			0.5 & 400
			& .098/.052/.010 & .122/.060/.013 & .099/.050/.011
			& .160/.091/.027 & .150/.091/.024 \\
			\addlinespace[0.35em]
			0.8 & 50
			& .097/.050/.011 & .170/.093/.020 & .113/.060/.015
			& .203/.118/.036 & .166/.104/.034 \\
			0.8 & 100
			& .100/.050/.010 & .150/.081/.018 & .106/.055/.013
			& .181/.105/.029 & .158/.097/.031 \\
			0.8 & 150
			& .091/.045/.010 & .137/.066/.012 & .105/.051/.011
			& .167/.095/.025 & .147/.089/.025 \\
			0.8 & 200
			& .094/.056/.010 & .133/.065/.015 & .098/.046/.012
			& .166/.100/.030 & .149/.087/.029 \\
			0.8 & 250
			& .094/.047/.009 & .128/.063/.013 & .102/.052/.010
			& .174/.100/.027 & .153/.089/.026 \\
			0.8 & 300
			& .094/.047/.009 & .134/.071/.016 & .098/.050/.008
			& .169/.096/.029 & .147/.085/.023 \\
			0.8 & 350
			& .090/.044/.009 & .122/.060/.011 & .094/.049/.010
			& .174/.101/.030 & .147/.089/.029 \\
			0.8 & 400
			& .091/.046/.009 & .130/.070/.012 & .098/.048/.010
			& .167/.097/.027 & .149/.085/.023 \\
			\bottomrule
		\end{tabular*}

		\vspace{0.35em}
		\begin{tablenotes}[flushleft]
			\scriptsize
			\item \emph{Notes.} The table reports empirical rejection frequencies under the null with \(B=5{,}000\) replications. Each entry is reported as \(10\%/5\%/1\%\), corresponding to the nominal significance level \(\alpha\). \(Q_{2}\) denotes the proposed vector-valued SN-CPA statistic, \(Q_{1}\) denotes the proposed scalar SN-CPA statistic, \(T_{\mathrm{SN}}\) denotes the unconditional self-normalized DM statistic, \(T_{\mathrm{GW}}\) denotes the HAC-based CPA statistic, and \(T_{\mathrm{DM}}\) denotes the HAC-based DM statistic.
		\end{tablenotes}
	\end{threeparttable}
\end{table}

\begin{center}
	\includegraphics[width=1.0\textwidth]{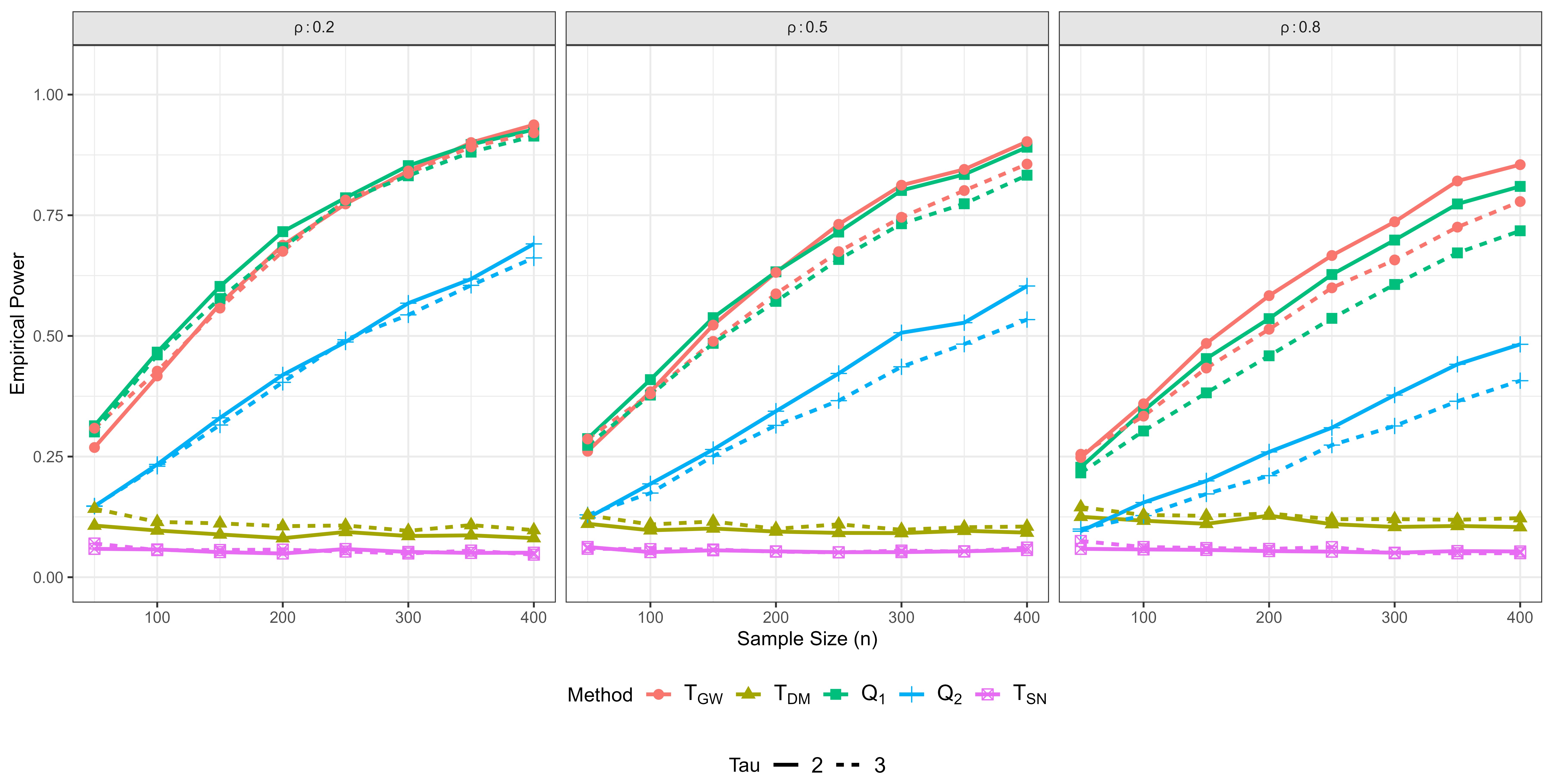}
	\captionof{figure}{Empirical power under DGP 1 for multistep forecasts.} 
	\label{fig:dgp1-power}
	\begin{minipage}{0.92\textwidth}
		\footnotesize
		Notes: The figure plots empirical rejection frequencies at the 5\% nominal level under the alternative with \(\delta=0.2\) and \(B=5{,}000\) replications. The panels correspond to different values of the persistence parameter \(\rho\). Solid and dashed lines correspond to \(\tau=2\) and \(\tau=3\), respectively.
	\end{minipage}
\end{center}

\section{Conclusion}\label{sec:conclusion}
In this paper, we propose a unified self-normalized framework for testing multistep conditional predictive ability. 
By employing functionals of the CUSUM process of the transformed loss differential, our approach avoids direct long-run covariance estimation and the associated ad hoc bandwidth selections required by conventional HAC methods. 
We establish asymptotic theory for both scalar and vector-valued test functions, deriving pivotal limiting null distributions, proving consistency against fixed alternatives, and obtaining noncentral limits under local alternatives. 
Monte Carlo simulations show that the proposed SN-CPA tests substantially alleviate the pronounced, horizon-dependent size distortions that affect traditional HAC-based methods, while maintaining robust empirical power.

\bibliographystyle{plainnat}
\bibliography{test}

\makeatletter
\@setaddresses
\global\let\addresses\@empty
\makeatother

\clearpage
\appendix
\setcounter{page}{1}
\renewcommand{\thepage}{S\arabic{page}}

\begin{center}
{\Large\bfseries Supplementary Material}
\end{center}

\section{Primitive Conditions for the High-Level Assumptions}

This section gives a brief illustration of how the high-level assumptions
imposed on the transformed loss-differential process $Z_{m,t+\tau}=h_t\Delta L_{m,t+\tau}$ can be supported by primitive conditions on the underlying data process, the forecasting methods, the loss function, and the test function. 
The discussion is only meant to clarify the link with primitive assumptions; the main results are stated under the high-level conditions used in the text.

\subsection{Mixing properties}

Suppose that the underlying process $\{W_t\}$ is $\alpha$-mixing with mixing coefficients of a suitable size.
Under the finite rolling window scheme, the forecasts $\widehat f_t$ and $\widehat g_t$ are measurable functions of a finite block of observations. 
For a fixed forecast horizon $\tau$, the loss differential $\Delta L_{m,t+\tau}= L(Y_{t+\tau},\widehat f_t)-L(Y_{t+\tau},\widehat g_t)$ is therefore a measurable function of a finite number of leads and lags of $\{W_t\}$. 
If, in addition, the test function $h_t$ is a finite-memory $\mathscr F_t$-measurable function, then $Z_{m,t+\tau}=h_t\Delta L_{m,t+\tau}$ is also a measurable function of a finite block of $\{W_t\}$. 
Standard mixing inheritance results for finite-block measurable transformations, such as Lemma 2.1 in \citet{WD84}, then imply that $\{Z_{m,t+\tau}\}$ is mixing with coefficients controlled by those of $\{W_t\}$, up to a finite lag shift.

\subsection{Moment conditions}

The moment conditions imposed on $Z_{m,t+\tau}$ can be verified from primitive moment bounds on the loss differential and the test function. 
For example, for each component $i$, suppose that, for some $\eta>0$,
\begin{equation*}
	\sup_t \mathbb E |h_{i,t}|^{2(2+\eta)}<\infty
\quad\text{and}\quad
\sup_t \mathbb E |\Delta L_{m,t+\tau}|^{2(2+\eta)}<\infty.
\end{equation*}
Then, by the Cauchy-Schwarz inequality,
\begin{equation*}
\mathbb E |Z_{i, m,t+\tau}|^{2+\eta} = \mathbb E |h_{i,t}\Delta L_{m,t+\tau}|^{2+\eta} \leq  \left(\mathbb E |h_{i,t}|^{2(2+\eta)}\right)^{1/2}\left(\mathbb E |\Delta L_{m,t+\tau}|^{2(2+\eta)}\right)^{1/2} <\infty .
\end{equation*}
Thus the required moment boundedness for the components of $Z_{m,t+\tau}$ follows.

\subsection{Partial sum regularity}

The preceding arguments verify the basic measurability, dependence, and moment requirements. 
The FCLT conditions used in the main text additionally require the corresponding long-run variance or covariance stability. 
These stability conditions are therefore kept as high-level assumptions in the paper, in line with the self-normalized limit theory, which depends on the entire partial sum process of $Z_{m,t+\tau}$.

\section{One-step Conditional Predictive Ability Test}
\label{sec:scalar-theory}
For the one-step case $\tau=1$, we also provide the corresponding test statistics and asymptotic theory. 

When the forecast horizon $\tau=1$, the null hypothesis becomes $H_0 \colon \E [ \Delta L_{m, t+1} \mid \mathscr F_{t} ]=0$. 
Since the test function $h_t$ is chosen to be $\mathscr{F}_t$-measurable, we have 
\begin{equation*}
    \E[Z_{m,t+1} \mid \mathscr{F}_t] = \E[h_t \Delta L_{m,t+1} \mid \mathscr{F}_t] = h_t \E[\Delta L_{m,t+1} \mid \mathscr{F}_t] = 0.
\end{equation*}
Thus, $\{Z_{m,t+1}\}$ is a martingale difference sequence (MDS) with respect to the filtration $\{\mathscr F_{t+1}\}_{t\geq m}$, under the null hypothesis, fulfilling the core structural assumption for the one-step self-normalized test.

\subsection{Scalar one-step case: $q=1, \tau=1$}
\label{subsec:scalar-one-theory}
 Under $H_0$, $\{Z_{m,t+1}\}$ is an MDS, so the long-run variance reduces to the contemporaneous variance. Although this variance can be consistently estimated, the adjusted-range statistic provides the basic self-normalized construction used as a benchmark for the remaining cases.

Let $n=T-m$ and define the CUSUM process
\begin{equation*}
	T_{m,n}^{h}(r) =\frac{1}{\sqrt{n}} \sum_{t=m}^{m+\lfloor nr\rfloor-1} \bigl(Z_{m, t+1}-\overline Z_{m,n}\bigr), \qquad 0\le r\le 1,
\end{equation*}
where $\overline Z_{m,n} =n^{-1}\sum_{t=m}^{T-1} Z_{m,t+1}$. The adjusted range is $ R^h_{m,n}= \sup_{0\leq r\leq 1}T^h_{m,n}(r) - \inf_{0\leq r\leq 1}T^h_{m,n}(r)$.
The scalar one-step self-normalized statistic is
\begin{equation*}
        Q^h_{m,n} = \frac{n\overline Z_{m,n}^{\,2}}{(R^h_{m,n})^2}.
\end{equation*}

\begin{assumption} \label{ass:scalar-one}
(1) There exists $\eta>0$ such that $\sup_t \E|Z_{m,t+1}|^{2+\eta}<\infty$.
(2) There exists $\sigma^2\in (0,\infty)$ such that, for every $r\in[0,1]$, $n^{-1} \sum_{t=m}^{m+\lfloor nr\rfloor-1}\E [Z_{m,t+1}^2 \mid \mathscr F_t ]
	  \overset{p}{\longrightarrow} r\sigma^2$.
\end{assumption}

Assumption~\ref{ass:scalar-one} is a high-level condition for $Z_{m,t+1}=h_t\Delta L_{m,t+1}$, ensuring the regularity needed for the partial sum FCLT; see \cite{GW06}, \cite{Brown71}, and \cite[Chapter~3]{HallHeyde80}.

\begin{theorem}
\label{thm:scalar-one-null}
Suppose Assumption~\ref{ass:scalar-one} holds. Then, under \(H_0\),
$ n^{-1/2} \sum_{t=m}^{m+\lfloor n \cdot \rfloor-1} Z_{m,t+1}\Rightarrow \sigma B(\cdot)$ in $D([0,1],\mathbb{R})$.
Consequently, $T^h_{m,n}(\cdot)\Rightarrow \sigma\mathbb B(\cdot)$ in $D([0,1],\mathbb{R})$, and
\begin{equation*}
        Q^h_{m,n} \Rightarrow \frac{B(1)^2}{\left(   \sup_{0\leq r\leq 1}\mathbb B(r) -\inf_{0\leq r\leq 1}\mathbb B(r)\right)^2}.      
\end{equation*}
\end{theorem}

Unlike the one-step Wald statistic of \cite{GW06},  the adjusted-range statistic eliminates the unknown scale through the CUSUM process itself.

\begin{theorem}
\label{thm:scalar-one-alt}
Suppose the statistic \(Q^h_{m,n}\) is defined as above.
\begin{enumerate}
\item[(i)] Under \(H^{(1,1),F}_{A,h}\), if Assumption~\ref{ass:scalar-one} holds with $Z_{m,t+1}$ replaced by $V_{m,t+1}^{F}$, then $Q^h_{m,n}\xrightarrow{p}\infty$.

\item[(ii)] Under $H^{(1,1),L}_{A,h}$, if Assumption~\ref{ass:scalar-one} holds with \(Z_{m,t+1}\) replaced by $V_{m,t+1}^{L}$, then, with $J=\mu_L/\sigma$,
\begin{equation*}
	 Q^h_{m,n} \Rightarrow \frac{(B(1)+J)^2} {\left( \sup_{0\leq r\leq 1}\mathbb B(r)-
        \inf_{0\leq r\leq 1}\mathbb B(r) \right)^2}.
\end{equation*}
\end{enumerate}
\end{theorem}

\subsection{Multivariate one-step case: $q>1,\tau=1$} \label{subsec:vector-one-theory}

For $q>1$ and $\tau=1$, $\{Z_{m,t+1}\}$ is a vector MDS under $H_0$. Serial autocovariances vanish, but the contemporaneous covariance matrix is generally non-diagonal. We use an LDL transformation for decorrelation, and then apply adjusted-range normalization componentwise.

Let $\widehat\Sigma_n= n^{-1}\sum_{t=m}^{T-1}Z_{m,t+1}Z_{m,t+1}^T$ and let $ \widehat\Sigma_n= \widehat D_n\widehat\Lambda_n\widehat D_n^T$ be its LDL decomposition, where $\widehat D_n$  is unit lower triangular and $\widehat\Lambda_n$ is diagonal with positive diagonal entries.
Define $ \widehat u_{m, t+1} =\widehat D_n^{-1}Z_{m,t+1}$ and  $ \overline{\widehat u}_{m, n}= n^{-1}\sum_{t=m}^{T-1} \widehat u_{m,t+1}$. 
Write $\widehat u_{m,t+1} =(\widehat u_{1,m,t+1},\ldots,\widehat u_{q,m,t+1})^T$  and 
    $\overline{\widehat u}_{m,n} = (\overline{\widehat u}_{1,m,n},\ldots \overline{\widehat u}_{q,m,n} )^T$.
For each $j=1,\ldots,q$, define $\widehat T_{j,m,n}^{h,q}(r) = n^{-1/2} \sum_{t=m}^{m+\lfloor nr \rfloor-1} \bigl( \widehat u_{j, m, t+1}-\overline{\widehat u}_{j,m,n}\bigr)$, $0\leq r \leq 1$ and  $ \widehat R^{h,q}_{j,m,n} = \sup_{0\leq r\leq 1}\widehat T^{h,q}_{j,m,n}(r)- \inf_{0\leq r\leq 1}\widehat T^{h,q}_{j,m,n}(r)$.
The multivariate one-step statistic is
\begin{equation*}
        Q^{h,q}_{m,n} =n\overline{\widehat u}_{m,n}^T \operatorname{diag} \left\{(\widehat R^{h,q}_{1,m,n})^{-2},\ldots,(\widehat R^{h,q}_{q,m,n})^{-2}\right\}\overline{\widehat u}_{m,n}.
\end{equation*}

\begin{assumption} \label{ass:vector-one}
(1) There exists a positive definite matrix $\Sigma$ such that $\Sigma=D\Lambda D^T$ 
where $D$ is unit lower triangular and $\Lambda=\mathrm{diag}\{\lambda_1, \ldots, \lambda_q \}$ with $\lambda_j>0$ for each $j=1,\ldots, q$. Moreover, $\widehat\Sigma_n\xrightarrow{p}\Sigma$, $\widehat D_n\xrightarrow{p}D$ and  $\widehat\Lambda_n\xrightarrow{p}\Lambda$.

(2) Let $u_{m,t+1}=D^{-1}Z_{m,t+1}$. There exists $\eta>0$ such that  $\sup_t \E \|u_{m,t+1}\|^{2+\eta}<\infty$, and, for every $r\in[0,1]$, $n^{-1} \sum_{t=m}^{m+\lfloor nr \rfloor-1} \E[ u_{m,t+1} u_{m,t+1}^T \mid \mathscr F_t] \overset{p}{\longrightarrow} r\Lambda$.
\end{assumption}

Assumption~\ref{ass:vector-one} is the multivariate martingale-FCLT condition after LDL decorrelation; see \cite{Brown71} and \cite[Chapter 3]{HallHeyde80}. It is the one-step vector analogue of the covariance regularity condition in \cite{GW06}.

\begin{theorem}
\label{thm:vector-one-null}
Suppose Assumption~\ref{ass:vector-one} holds. Then, under $H_0$, $
       n^{-1/2} \sum_{t=m}^{m+\lfloor n \cdot \rfloor-1} \widehat u_{m,t+1}  \Rightarrow\Lambda^{1/2}B^q(\cdot)$ in $D([0,1],\mathbb{R}^q)$.
Consequently, for each \(j=1,\ldots,q\), $ \widehat T^{h,q}_{j,m,n}(\cdot)\Rightarrow \sqrt{\lambda_j}\,\mathbb B^q_j(\cdot)$ in $D([0,1],\mathbb{R})$, and 
\begin{equation*}
        Q^{h,q}_{m,n} \Rightarrow \sum_{j=1}^q  \frac{B^q_j(1)^2} {\left( \sup_{0\leq r\leq 1}\mathbb B^q_j(r)- \inf_{0\leq r\leq 1}\mathbb B^q_j(r)  \right)^2}.       
\end{equation*}
\end{theorem}

This is the multivariate one-step analogue of the scalar result. 
Unlike the multivariate Wald statistic in \cite{GW06}, this statistic does not use the inverse of an estimated covariance matrix for scale normalization.
The LDL step is used only for decorrelation, while adjusted ranges remove the coordinatewise scales.

\begin{assumption}
\label{ass:vector-one-fixed}
Under \(H^{(q,1),F}_{A,h}\), let $V_{m,t+1}^F=Z_{m,t+1}-\mu_F$. The partial sum process
$ \left\{ n^{-1/2} \sum_{t=m}^{m+\lfloor n \cdot \rfloor-1}V_{m,t+1}^F \right\}$ is tight in $\ell^\infty([0,1],\mathbb R^q)$. Moreover, $\widehat\Sigma_n\xrightarrow{p}\Sigma_F$ where $\Sigma_F$ is positive definite.
  Let $\Sigma_F=D_F\Lambda_FD_F^T$ be its LDL decomposition, where \(D_F\) is unit lower triangular and $\Lambda_F$ is diagonal with strictly positive diagonal entries.
  \end{assumption}

\begin{assumption} \label{ass:vector-one-local}
Under \(H^{(q,1),L}_{A,h}\), Assumption~\ref{ass:vector-one} holds with $ W_{m,t+1}
        = D^{-1}Z_{m,t+1}- n^{-1/2}D^{-1}\mu_L$ in place of $u_{m,t+1}$.
\end{assumption}

The fixed and local alternatives are treated separately because a fixed mean shift may change the probability limit of $\widehat\Sigma_n$, whereas an $n^{-1/2}$-local shift does not affect its first-order limit.

\begin{theorem} \label{thm:vector-one-alt}
Suppose the statistic $Q^{h,q}_{m,n}$ is defined as above.
\begin{enumerate}
\item[(i)] Under \(H^{(q,1),F}_{A,h}\), if
Assumption~\ref{ass:vector-one-fixed} holds, then $ Q^{h,q}_{m,n}\xrightarrow{p}\infty$.
\item[(ii)] Under $H^{(q,1),L}_{A,h}$, if Assumption~\ref{ass:vector-one-local} holds, then, with $J=\Lambda^{-1/2}D^{-1}\mu_L=(J_1,\ldots,J_q)^T$, we have
\begin{equation*}
	Q^{h,q}_{m,n} \Rightarrow  \sum_{j=1}^q  \frac{(B^q_j(1)+J_j)^2} {\left( \sup_{0\leq r\leq 1}\mathbb B^q_j(r)-  \inf_{0\leq r\leq 1}\mathbb B^q_j(r)\right)^2}.
\end{equation*} 
\end{enumerate}
\end{theorem}

\clearpage
\section{Simulation results for DGP2}

Tables~\ref{tab:dgp2-size-tau2} and~\ref{tab:dgp2-size-tau3} report empirical
size under DGP~2. The proposed vector-valued  statistic \(Q_2\) remains close to the
nominal levels across the multistep designs, whereas the HAC-based statistics
tends to over-reject, especially for \(\tau=3\).

\begin{table}[H]
	\centering
	\begin{threeparttable}
		\caption{Empirical size under DGP 2 for \(\tau=2\)}
		\label{tab:dgp2-size-tau2}

		\fontsize{7.5}{8.3}\selectfont
		\setlength{\tabcolsep}{5.0pt}
		\renewcommand{\arraystretch}{1.02}

		\begin{tabular*}{\textwidth}{@{\extracolsep{\fill}}cc ccccc@{}}
			\toprule
			\(p\) & \(n\)
			& \(Q_2\) & \(Q_1\) & \(T_{\mathrm{SN}}\)
			& \(T_{\mathrm{GW}}\) & \(T_{\mathrm{DM}}\) \\
			\midrule
			0.2 & 50
			& .098/.049/.010 & .179/.091/.016 & .111/.056/.013
			& .162/.096/.030 & .170/.097/.034 \\
			0.2 & 100
			& .102/.053/.010 & .161/.082/.014 & .104/.056/.012
			& .148/.085/.023 & .164/.097/.031 \\
			0.2 & 150
			& .103/.053/.011 & .151/.078/.016 & .107/.057/.011
			& .147/.080/.021 & .157/.090/.025 \\
			0.2 & 200
			& .099/.047/.009 & .145/.076/.015 & .102/.051/.009
			& .142/.079/.020 & .154/.090/.025 \\
			0.2 & 250
			& .102/.053/.011 & .141/.081/.018 & .103/.055/.011
			& .147/.082/.021 & .154/.090/.025 \\
			0.2 & 300
			& .099/.055/.014 & .131/.066/.014 & .103/.052/.011
			& .140/.079/.021 & .147/.086/.028 \\
			0.2 & 350
			& .093/.050/.010 & .139/.069/.013 & .095/.050/.011
			& .136/.076/.018 & .146/.088/.022 \\
			0.2 & 400
			& .095/.050/.011 & .127/.060/.012 & .104/.051/.011
			& .132/.075/.021 & .150/.087/.026 \\
			\addlinespace[0.35em]

			0.5 & 50
			& .103/.057/.014 & .174/.097/.028 & .110/.056/.015
			& .160/.098/.035 & .167/.099/.035 \\
			0.5 & 100
			& .095/.048/.009 & .154/.085/.021 & .105/.053/.009
			& .154/.087/.025 & .157/.088/.030 \\
			0.5 & 150
			& .098/.055/.012 & .148/.080/.017 & .100/.053/.010
			& .149/.085/.024 & .155/.090/.024 \\
			0.5 & 200
			& .107/.054/.010 & .137/.074/.016 & .104/.053/.011
			& .135/.077/.026 & .155/.091/.025 \\
			0.5 & 250
			& .097/.050/.011 & .130/.066/.015 & .100/.050/.010
			& .138/.076/.022 & .148/.083/.025 \\
			0.5 & 300
			& .094/.051/.010 & .123/.065/.015 & .097/.049/.011
			& .139/.077/.020 & .144/.082/.025 \\
			0.5 & 350
			& .098/.047/.009 & .119/.061/.013 & .099/.052/.011
			& .130/.075/.019 & .147/.085/.023 \\
			0.5 & 400
			& .094/.055/.007 & .127/.065/.015 & .103/.050/.011
			& .137/.079/.019 & .151/.090/.027 \\
			\addlinespace[0.35em]

			0.8 & 50
			& .098/.047/.007 & .177/.097/.027 & .106/.058/.012
			& .164/.098/.031 & .161/.097/.031 \\
			0.8 & 100
			& .096/.050/.008 & .149/.082/.020 & .108/.054/.012
			& .147/.082/.025 & .156/.094/.031 \\
			0.8 & 150
			& .099/.045/.009 & .140/.074/.016 & .107/.050/.009
			& .142/.082/.022 & .153/.086/.026 \\
			0.8 & 200
			& .099/.052/.008 & .136/.070/.013 & .105/.055/.012
			& .143/.078/.019 & .147/.082/.023 \\
			0.8 & 250
			& .100/.054/.011 & .127/.070/.016 & .101/.054/.012
			& .150/.081/.022 & .151/.092/.028 \\
			0.8 & 300
			& .103/.054/.013 & .128/.066/.014 & .108/.053/.012
			& .143/.077/.021 & .152/.092/.025 \\
			0.8 & 350
			& .100/.050/.009 & .118/.061/.013 & .099/.053/.010
			& .133/.076/.017 & .149/.086/.023 \\
			0.8 & 400
			& .103/.056/.010 & .130/.066/.016 & .106/.053/.012
			& .143/.077/.023 & .153/.094/.029 \\
			\bottomrule
		\end{tabular*}

		\vspace{0.35em}
		\begin{tablenotes}[flushleft]
			\scriptsize
			\item \emph{Notes.} The table reports empirical rejection frequencies under the null with \(B=5{,}000\) replications. Each entry is reported as \(10\%/5\%/1\%\), corresponding to the nominal significance level \(\alpha\). \(Q_2\) denotes the proposed vector-valued SN-CPA statistic, \(Q_1\) denotes the proposed scalar SN-CPA statistic, \(T_{\mathrm{SN}}\) denotes the unconditional self-normalized DM statistic, \(T_{\mathrm{GW}}\) denotes the HAC-based CPA statistic, and \(T_{\mathrm{DM}}\) denotes the HAC-based DM statistic.
		\end{tablenotes}
	\end{threeparttable}
\end{table}

\begin{table}[H]
	\centering
	\begin{threeparttable}
		\caption{Empirical size under DGP 2 for \(\tau=3\)}
		\label{tab:dgp2-size-tau3}

		\fontsize{7.5}{8.3}\selectfont
		\setlength{\tabcolsep}{5.0pt}
		\renewcommand{\arraystretch}{1.02}

		\begin{tabular*}{\textwidth}{@{\extracolsep{\fill}}cc ccccc@{}}
			\toprule
			\(p\) & \(n\)
			& \(Q_2\) & \(Q_1\) & \(T_{\mathrm{SN}}\)
			& \(T_{\mathrm{GW}}\) & \(T_{\mathrm{DM}}\) \\
			\midrule
			0.2 & 50
			& .112/.057/.011 & .197/.099/.019 & .126/.066/.016
			& .206/.125/.048 & .207/.131/.054 \\
			0.2 & 100
			& .103/.048/.010 & .166/.085/.016 & .104/.057/.013
			& .175/.102/.031 & .179/.115/.038 \\
			0.2 & 150
			& .097/.053/.011 & .158/.080/.018 & .106/.055/.014
			& .162/.095/.029 & .167/.102/.035 \\
			0.2 & 200
			& .097/.048/.009 & .147/.079/.015 & .096/.047/.010
			& .146/.086/.023 & .163/.098/.032 \\
			0.2 & 250
			& .107/.055/.011 & .145/.079/.020 & .109/.057/.012
			& .167/.094/.028 & .175/.103/.031 \\
			0.2 & 300
			& .097/.050/.010 & .133/.072/.015 & .098/.047/.010
			& .156/.088/.023 & .174/.104/.031 \\
			0.2 & 350
			& .102/.051/.010 & .137/.067/.016 & .100/.046/.009
			& .155/.091/.024 & .168/.103/.031 \\
			0.2 & 400
			& .096/.048/.009 & .137/.070/.014 & .109/.054/.009
			& .146/.083/.025 & .166/.099/.033 \\
			\addlinespace[0.35em]

			0.5 & 50
			& .115/.064/.013 & .197/.118/.034 & .119/.066/.017
			& .205/.130/.050 & .199/.126/.048 \\
			0.5 & 100
			& .103/.055/.010 & .170/.099/.022 & .108/.053/.012
			& .182/.108/.035 & .183/.113/.043 \\
			0.5 & 150
			& .103/.055/.010 & .153/.087/.020 & .109/.060/.013
			& .175/.105/.033 & .184/.113/.038 \\
			0.5 & 200
			& .106/.056/.012 & .150/.086/.023 & .115/.060/.011
			& .162/.099/.032 & .178/.107/.038 \\
			0.5 & 250
			& .103/.051/.010 & .140/.075/.019 & .098/.054/.012
			& .153/.088/.028 & .171/.106/.032 \\
			0.5 & 300
			& .108/.057/.010 & .135/.070/.017 & .107/.059/.014
			& .161/.092/.027 & .181/.107/.037 \\
			0.5 & 350
			& .098/.049/.010 & .139/.071/.016 & .102/.054/.009
			& .157/.095/.026 & .176/.101/.035 \\
			0.5 & 400
			& .103/.049/.009 & .133/.067/.015 & .102/.051/.008
			& .154/.093/.025 & .173/.106/.030 \\
			\addlinespace[0.35em]

			0.8 & 50
			& .119/.063/.013 & .211/.130/.041 & .129/.073/.018
			& .210/.136/.050 & .212/.136/.055 \\
			0.8 & 100
			& .104/.056/.011 & .175/.095/.026 & .113/.061/.016
			& .169/.105/.036 & .189/.112/.040 \\
			0.8 & 150
			& .105/.058/.012 & .149/.083/.023 & .105/.054/.014
			& .163/.098/.027 & .171/.107/.038 \\
			0.8 & 200
			& .107/.055/.010 & .154/.084/.022 & .116/.057/.012
			& .168/.099/.033 & .186/.112/.041 \\
			0.8 & 250
			& .098/.048/.010 & .135/.072/.016 & .101/.054/.011
			& .159/.091/.029 & .168/.102/.034 \\
			0.8 & 300
			& .100/.056/.012 & .136/.070/.017 & .108/.055/.012
			& .169/.096/.029 & .182/.115/.036 \\
			0.8 & 350
			& .097/.050/.010 & .133/.072/.016 & .107/.053/.011
			& .158/.090/.025 & .171/.104/.032 \\
			0.8 & 400
			& .101/.054/.009 & .127/.070/.012 & .103/.054/.011
			& .156/.091/.024 & .171/.106/.032 \\
			\bottomrule
		\end{tabular*}

		\vspace{0.35em}
		\begin{tablenotes}[flushleft]
			\scriptsize
			\item \emph{Notes.} The table reports empirical rejection frequencies under the null with \(B=5{,}000\) replications. Each entry is reported as \(10\%/5\%/1\%\), corresponding to the nominal significance level \(\alpha\). \(Q_{2}\) denotes the proposed vector-valued SN-CPA statistic, \(Q_{1}\) denotes the proposed scalar SN-CPA statistic, \(T_{\mathrm{SN}}\) denotes the unconditional self-normalized DM statistic, \(T_{\mathrm{GW}}\) denotes the HAC-based CPA statistic, and \(T_{\mathrm{DM}}\) denotes the HAC-based DM statistic.
		\end{tablenotes}
	\end{threeparttable}
\end{table}
\begin{figure}[t]
\centering
\includegraphics[width=1\textwidth]{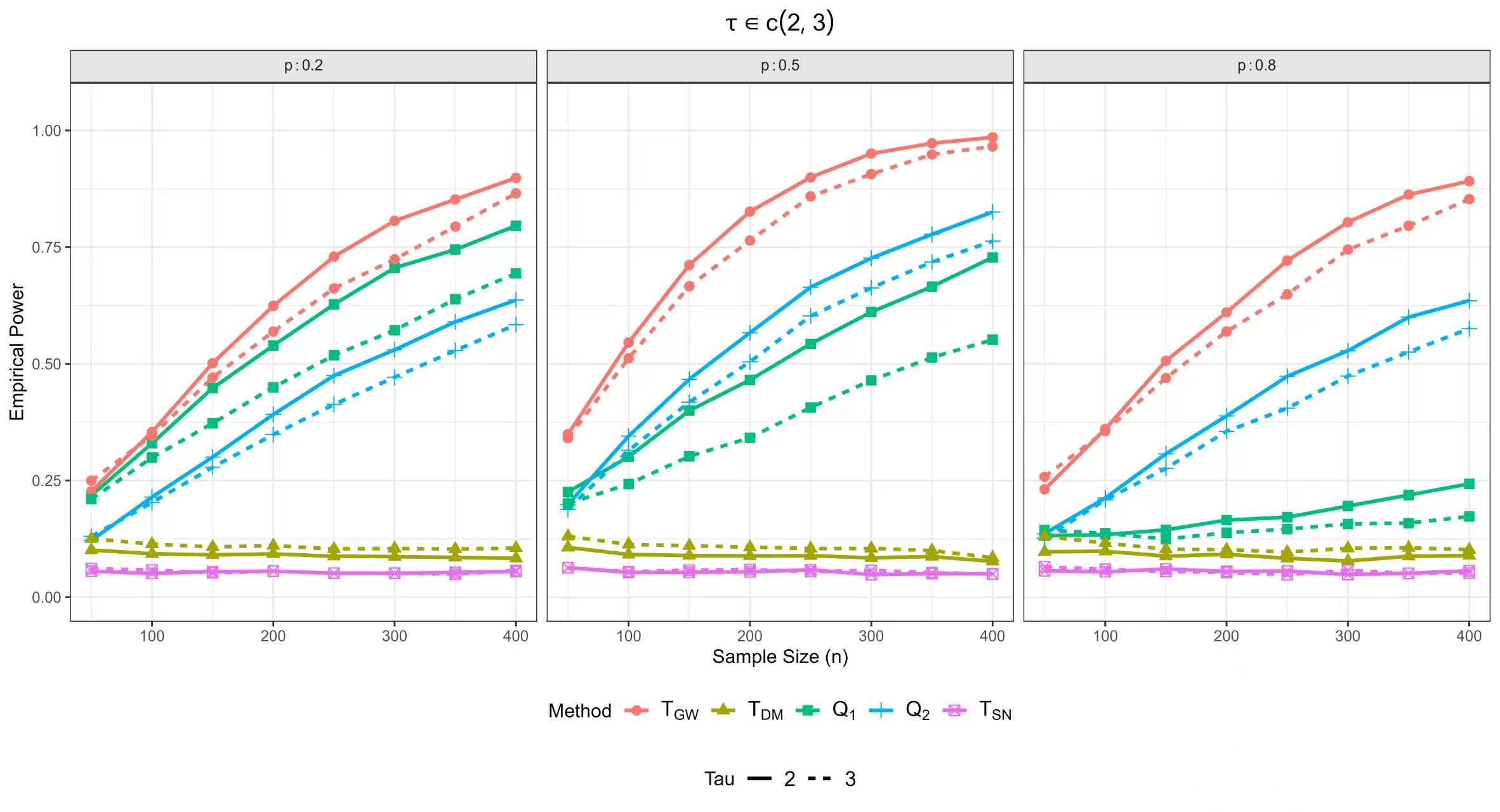}
\caption{Empirical power under DGP 2 for multistep forecasts.}
\label{fig:dgp2-power}
\begin{minipage}{0.92\textwidth}
\footnotesize
Notes: The figure plots empirical rejection frequencies at the 5\% nominal level under the alternative with \(B=5{,}000\) replications. The panels correspond to different values of \(p\). Solid and dashed lines correspond to \(\tau=2\) and \(\tau=3\), respectively.
\end{minipage}
\end{figure}

Figure~\ref{fig:dgp2-power} shows that the power of the proposed \(Q_2\)
statistic increases with the evaluation sample size ($d=0.5$). The unconditional
self-normalized DM statistic has limited power because the alternatives are
conditional rather than unconditional.

\end{document}